\documentclass[11pt]{amsart}
\usepackage{amssymb}

\usepackage{palatino}
\input amssym.def

\usepackage{amsmath, amsfonts}
\usepackage{amssymb}
\usepackage{amscd}
\usepackage[mathscr]{eucal}
\usepackage{palatino}
\newfont{\cyrr}{wncyr10}

\newcommand{\N}{{\mathbb N}}
\newcommand{\Z}{{\mathbb Z}}

\newcommand{\R}{{\mathbb R}}

\newcommand{\SL}{{\rm SL}}

\newtheorem{thm}{Theorem}

\newtheorem{lem}[thm]{Lemma}
\newtheorem{cor}{Corollary}

\newtheorem{rmk}{Remark}[section]

\begin{document}

\title[Omega result]{On Hecke eigenvalues of Siegel modular forms
in the Maass space}

\author{Sanoli Gun, Biplab Paul and Jyoti Sengupta}

\address[Sanoli Gun and Biplab Paul]   
{Institute of Mathematical Sciences, HBNI,
C.I.T Campus, Taramani, 
Chennai  600 113, 
India.}
\email{sanoli@imsc.res.in}
\email{biplabpaul@imsc.res.in}

\address[J. Sengupta]   
{School of Mathematics,
Tata institute of Fundamental research, 
Homi Bhabha Road, Mumbai 400 005, 
India.}
\email{sengupta@math.tifr.res.in}



\begin{abstract} 
In this article, we prove an omega-result for the 
Hecke eigenvalues $\lambda_F(n)$ of Maass forms $F$
which are Hecke eigenforms in the space of Siegel 
modular forms of weight $k$, genus two for
the Siegel modular group $Sp_2(\Z)$. 
In particular, we prove
$$
\lambda_F(n)= \Omega(n^{k-1}\text{exp}
(c \frac{\sqrt{\log n}}{\log\log n})),
$$
when $c>0$ is an absolute constant. 
This improves the earlier result 
$$
\lambda_F(n)= \Omega(n^{k-1} 
(\frac{\sqrt{\log n}}{\log\log n}))
$$
of Das and the third author. We also show that
for any $n \ge 3$, one has
$$
\lambda_F(n) \leq n^{k-1}\text{exp}
\left(c_1\sqrt{\frac{\log n}{\log\log n}}\right),
$$
where $c_1>0$ is an absolute constant.
This improves an earlier result of Pitale and Schmidt.
Further, we investigate the limit
points of the sequence 
$\{\frac{\lambda_F(n)}{n^{k-1}}\}_{n \in \N}$
and show that it has infinitely many 
limit points.  Finally, we show that $\lambda_F(n) >0$
for all $n$, a result earlier proved by
Breulmann by a different technique.
\end{abstract}

\maketitle

\section{Introduction}

For $g \ge 1$, let $\Gamma_g := Sp_g(\Z)$ be the 
Siegel modular group of genus $g$ and 
$S_k^g$ be the space of cuspidal Siegel modular forms 
of weight $k$ and genus $g$ for $\Gamma_g$.
One of the interesting problem in the theory of modular forms 
is to understand arithmetic nature of eigenvalues of the 
Hecke operators
acting on the space $S_k^g$. 
Let $f$ be a normalised Hecke eigenform of weight $k$ and genus 
$g=1$ with the Hecke eigenvalues $\lambda_f(n)$. By a
celebrated result of Deligne, one has 
\begin{equation*}
|\lambda_f(n)|\leq d(n)~n^{(k-1)/2},
\end{equation*} 
where $d(n)$ is the number of divisors of $n$.
One would like to know the optimality of the 
above result, i.e. an omega result for the sequence
$\{\lambda_f(n)/n^{(k-1)/2}\}_{n\in \N}$. 
In 1973, Rankin \cite{RR} showed that
$$
\limsup_{n} \frac{\lambda_f(n)}{n^{(k-1)/2}} = + \infty.
$$
In 1983, Ram Murty \cite{RM} showed that 
\begin{equation*}
|\lambda_f(n)|=\Omega\left(n^{(k-1)/2}\text{exp}
(c\log n/\log\log n)\right),
\end{equation*}
where $c>0$ is an absolute constant. 

It is natural to investigate similar questions for 
higher genus. In this direction, the generalised 
Ramanujan-Petersson 
conjecture \cite{AP} implies that for any prime $p$
and $\epsilon>0$, one has
\begin{equation}\label{eq1}
\lambda_F(p)\ll_{g,\epsilon}p^{gk/2-g(g+1)/4+\epsilon}.
\end{equation}
However, it is known that when $g=2$,
the elements of the Maass space in $S_k^2$ are
precisely the ones which fail to satisfy equation 
\eqref{eq1}. R. Weissauer \cite{RW} showed that
Ramanujan-Petersson conjecture is true when $F$ 
does not belong to the Maass space in $S_k^2$.

From now on, we concentrate on the space of Maass cusp 
forms denoted by $S_k^*$ in the space $S_k^2$. In a 
recent work \cite{DS},  Das along with the third author 
studied the question of omega result for Hecke 
eigenvalues of the Hecke operators acting on
$S_k^*$. But there is a considerably large gap 
between the known upper bound 
$$
\lambda_F(n) \ll_{\epsilon} n^{k-1 + \epsilon},
\phantom{m}\text{ for any } \epsilon >0
$$
and the known omega result. 

In this article, we investigate arithmetic
behaviour of Hecke eigenvalues
of Maas forms in $S_k^*$ and also study 
the limit points of the sequence
$\{\lambda_F(n)/n^{k-1} \}_{ n \in \N}$.
More precisely, we prove the following theorems.

\begin{thm}\label{th1}
Let $F\in S_k^*$ be a non-zero Hecke eigenform. Then 
there exists an absolute constant $c>0$ such that 
\begin{equation}
\lambda_F(n)=\Omega(n^{k-1}\rm{exp}
(c\frac{\sqrt{\log n}}{\log\log n}))
\end{equation}
\end{thm}

Our next theorem shows that the above omega result
is not too far from an upper bound one can derive. 
In particular, we have
\begin{thm}\label{th2}
Let $F\in S_k^*$ be a non-zero Hecke eigenform. Then 
there exists an absolute constant $c_1>0$ such that  
\begin{equation}
\lambda_F(n) ~\leq~ n^{k-1}\rm{exp}
\left(c_1\sqrt{\frac{\log n}{\log\log n}}\right)
\end{equation}
for all $n\in \N$ with $n\ge 3$.
\end{thm}

\begin{rmk}
Theorem \ref{th2} improves an earlier result of
Pitale and Schmidt (see page 101 of \cite{PS}).
\end{rmk}

We also have the following lower bound.
\begin{thm}\label{th3}
Let $F\in S_k^*$ be a non-zero Hecke eigenform. 
Then there exist absolute constants $c_2, c_3>0$ such that 
\begin{equation}
\lambda_F(n)\ge c_2 n^{k-1}\rm{exp}
\left(-c_3\sqrt{\frac{\log n}{\log\log n}}\right)
\end{equation}
for all $n\in \N$ with $n \ge 3$.
\end{thm}

As a corollary, we now get the following result of Breulmann \cite{SB}
whose proof is different from ours.
\begin{cor}\label{co}
If $F \in S_k^*$ is a non-zero Hecke eigenform with Hecke eigenvalues 
$\lambda_F(n)$, then $\lambda_F(n)>0$.  
\end{cor}

Since $\frac{\lambda_F(n)}{n^{k-1}}>0$, one can ask whether this 
result is optimal. Our next theorem shows that the 
answer is positive.
 
\begin{thm}\label{th4}
Let $F\in S_k^*$ be a non-zero Hecke eigenform. Then 
\begin{equation}
\liminf_{n}
\frac{\lambda_F(n)}{n^{k-1}}=0
\end{equation}
\end{thm}

Finally, we investigate the limit points of the
sequence $\{ \lambda_F(n)/n^{k-1}\}_{n \in \N}$. 
In this direction, we have the following result.
\begin{thm}\label{th5}
Infinitely many limit points of the sequence 
$\{ \frac{\lambda_F(n)}{n^{k-1}} \}_{ n \in \N}$ are greater
than $1$ and infinitely many of them are less than $1$.
\end{thm}
In order to prove our results, we rely on an
idea of Rankin \cite{RR} and some 
standard analytic techniques. We manage to avoid
the use of Sato-Tate conjecture which is now a 
theorem due to  Barnet-Lamb, Geraghty, Harris 
and Taylor \cite{BGHT}.

\section{Notation and preliminaries}\label{sec2}

Throughout the article, let $\mathcal{P}$ denote the set of all 
rational prime numbers. Also we use the notation $q:=e^{2\pi iz}$, 
where $z\in\mathcal{H}$, the complex upper half-plane. We say 
that $f(x)=\Omega(g(x))$ to indicate that $\limsup_{x
\rightarrow\infty}|f(x)|/g(x)>0$. Moreover, we shall write 
$f(x)\sim g(x)$ when $\lim_{x\rightarrow\infty}f(x)/g(x)=1$. 
A subset $A \subset \mathcal{P}$ 
is said to have the lower natural density $\alpha(A)$ if 
\begin{equation*}
\alpha(A)=\liminf_{x\rightarrow\infty}
\frac{\#\{p\leq x: p\in A\}}{\#\{p\leq x\}}
\end{equation*}

We now recall the following lemma which we 
will use to prove our results.
\begin{lem}(\cite{DS}, lemma 3.1)\label{lem6}
Let $f(z)=\sum_{n=1}^\infty a(n)q^n\in S_k^1$ be a 
normalised Hecke eigenform. Then there exists an absolute 
constant $0< \beta <2$ such that the set 
\begin{equation*}
\{p\in\mathcal{P}: a(p)> \beta\cdot p^{(k-1)/2}\}
\end{equation*}
has positive lower density.
\end{lem}
One can use the Sato-Tate conjecture to get the 
above result but proof of the lemma avoids that.

\section{Proof of theorem \ref{th1}}\label{sec3}

Let $F\in S_k^*$ be a nonzero Hecke eigenform with
eigenvalues $\lambda_F(n)$. Then there 
exists a normalised Hecke eigenform $f$ of weight $2k-2$ 
for the full modular group $\SL_2(\Z)$ such that $F$ is 
the Saito-Kurokawa lift of $f$. We know that for any prime $p$, 
one has (see \cite{SB} for details)
$$
\lambda_F(p)=p^{k-1}(1+\frac{1}{p}+\frac{a(p)}{p^{k-1}}).
$$
Note that by lemma \ref{lem6}, there exists an absolute 
constant $0< \beta <2$ such that the set
$$
A:=\{p: a(p)> \beta \cdot p^{k-3/2}\}
$$
has positive lower density. For any $x > 0$, let
$$
n_x :=\prod_{5 \leq p\leq x, \atop p\in A}p
$$
with the convention that an empty
product is $1$. Then for sufficiently large $x \in \R^+$, 
we have 
\begin{eqnarray*}
\frac{\lambda_F(n_x)}{n_x^{k-1}}
~~=~~  
\prod_{5 \leq p\leq x, \atop p\in A}
\left(1+\frac{1}{p}+\frac{a(p)}{p^{k-1}}\right)
~~\geq ~~
\prod_{5 \leq p\leq x,  \atop p\in A}
\left(1+\frac{a(p)}{p^{k-1}}\right)
& \geq & 
\text{exp}\left[\sum_{5 \leq p\leq x, \atop p\in A}
\log(1+\frac{\beta}{p^{1/2}})\right]\\
& \geq & 
\text{exp}\left(c_4\sum_{5 \leq p\leq x, \atop p\in A}
\frac{1}{p^{1/2}}\right),
\end{eqnarray*}
where $c_4 > 0$ is an absolute constant.
Since the set $A$ has positive lower density, by partial 
summation formula, it can be seen that 
$$
\sum_{5 \leq p\leq x, \atop p\in A}
\frac{1}{p^{1/2}}
~~\gg~~
\frac{\sqrt{x}}{\log x}~,
$$
where the implied constant is absolute. Further note that 
for any $x \in \R^+$, we have
\begin{equation*}
\log (n_x)
=
\sum_{5 \leq p\leq x, \atop p\in A}\log p
~\ll~
x
\end{equation*}
with the convention that an empty sum is zero.
Note that $\frac{\sqrt{x}}{{\log x}}$ is an increasing function
for $x \ge 8$. Thus for sufficiently large $x$, we have
\begin{eqnarray*}
\frac{\lambda_F(n_x)}{n_x^{k-1}}
~\geq~ 
\text{exp}\left(c_5\frac{\sqrt{x}}{\log x}\right)
~\geq~ 
\text{exp}\left(c\frac{\sqrt{\log n_x}}{\log\log n_x}
\right),
\end{eqnarray*}
where $c, c_5 > 0$ are absolute constants.
This shows that given any natural number $M$, there exists 
a  natural number $n$ with $n>M$ such that 
$$
\frac{\lambda_F(n)}{n^{k-1}}\geq \text{exp}
\left(c\frac{\sqrt{\log n}}{\log\log n}\right).
$$
This completes the proof of Theorem \ref{th1}.

\section{Proof of theorem \ref{th2}}\label{sec4}

In this section, we keep the notations as in
the section \ref{sec2} and section \ref{sec3}. It can be
deduced from \cite{SB}
that for all $m\in\N$ and any $p \in \mathcal{P}$, we have
\begin{equation*}
\frac{\lambda_F(p^m)}{p^{m(k-1)}}
~=~
1 + \frac{1}{p} ~+~(1+\frac{1}{p})
\sum_{\ell=1}^{m-1}\frac{a(p^{\ell})}{p^{\ell(k-1)}}
~+~ \frac{a(p^m)}{p^{m(k-1)}}
\end{equation*}
with the convention that an empty sum is
zero. For any $ |\lambda| < 1$, the series
$$
\sum_{n=2}^{\infty} (n+1)\lambda^{n}
~=~
\sum_{n=3}^{\infty} n  \lambda^{n-1}
~=~
\frac{3\lambda^2 - 2 \lambda^3}{(1 - \lambda)^2}.
$$ 
This can be seen by considering the power series
$$
f(Y) = \sum_{n \ge 3} Y^n = \frac{1}{1-Y} - 1 - Y - Y^2
$$
and noting that 
$$
f'(Y) ~=~
\frac{3Y^2 - 2 Y^3}{(1 - Y)^2},
$$
where $f'$ is the derivative of $f$.
For any $p \in \mathcal{P}$, let us
set
$$
\alpha_p 
:= \sum_{n=2}^{\infty}
\frac{n+1}{p^{{n}/2}}
~=~
\frac{3 p^{1/2} - 2}{p^{1/2}(p^{1/2} -1)^2}.
$$ 
By the work of Deligne, one knows that
$$
\frac{a(n)}{n^{k-3/2}}
~\leq~ d(n),
$$ 
where $d(n)$ denotes the number of divisors of $n$. 
This shows that for any $p \in \mathcal{P}$ 
and $m\in\N$ with $m \geq 2$, we have 
\begin{eqnarray*}
\frac{\lambda_F(p^m)}{p^{m(k-1)}} 
& \leq & 
1+\frac{1}{p}
 ~+~
(1+\frac{1}{p})\frac{2}{p^{1/2}}
~+~ (1+\frac{1}{p})\alpha_p.
\end{eqnarray*}
Note that $\alpha_p \asymp\frac{1}{p}$. 
Hence there exists an absolute constant $c_7>0$ such that 
\begin{eqnarray*}\label{eq9}
\frac{\lambda_F(p^m)}{p^{m(k-1)}} 
~\leq~  1+\frac{c_7}{p^{1/2}} 
\end{eqnarray*}
for all $m \in \N$.
Let $n \ge 3$ be an arbitrary natural number and let $t = \nu(n)$
be its number of distinct prime divisors. Then we can write  $n$ as 
$$
n= p_1^{m_1} \cdots p_t^{m_t}
$$ 
where $p_1 < \cdots < p_t$ and $m_i >0$ for $1 \le i \le t$. 
Thus we have 
\begin{equation*}
\frac{\lambda_F(n)}{n^{k-1}} 
~~\leq~
\prod_{1\leq i\leq t}(1+\frac{c_7}{p_i^{1/2}})
~=~ 
\text{exp}\left(\sum_{1\leq i\leq t}
\log(1+\frac{c_7}{p_i^{1/2}})\right) 
~\leq~  
\text{exp}\left(c_7\sum_{1\leq i\leq t}
\frac{1}{p_i^{1/2}}\right).
\end{equation*}
Here we have used the fact that
$\log(1+x) \le x$ for any $x >0$.
Since $i < p_i$, we have
\begin{eqnarray*}
\frac{\lambda_F(n)}{n^{k-1}} 
& \leq & 
\text{exp}\left(c_7\sum_{1\leq i\leq t}
\frac{1}{i^{1/2}}\right)\leq\text{exp}\left(c_8~ t^{1/2}
\right),
\end{eqnarray*}
where $c_8 >0$ is an absolute constant.
Note that $t = \nu(n) \ll \frac{\log n}{\log\log n}$
for $n \gg 1$(see \cite{GT}, page 83 for details).   
Thus for any $n \ge 3$, we have 
\begin{equation}\label{eq10}
\frac{\lambda_F(n)}{n^{k-1}}\leq 
\text{exp}\left(c_1\sqrt{\frac{\log n}{\log\log n}}\right),
\end{equation}
where $c_1>0$ is an absolute constant. 
This completes the proof of the theorem.

\section{Proof of theorem \ref{th3}}\label{sec5}

As earlier, we keep the notations as in the previous sections.
We know that for any $p \in \mathcal{P}$, one has
\begin{equation}\label{imp}
\frac{\lambda_F(p^m)}{p^{m(k-1)}}
~=~
1 + \frac{1}{p} ~+~(1+\frac{1}{p})
\sum_{\ell=1}^{m-1}\frac{a(p^{\ell})}{p^{\ell(k-1)}}
~+~ \frac{a(p^m)}{p^{m(k-1)}}
\end{equation}
with the convention that an empty sum is zero.
Proceeding as in section \ref{sec4}, 
for any $p \in \mathcal{P}$ and $m\in\N$
with $m \ge 2$, we see that
\begin{eqnarray*}
\frac{\lambda_F(p^m)}{p^{m(k-1)}} 
& \geq & 
1 + \frac{1}{p}
~+~
(1+\frac{1}{p})\frac{a(p)}{p^{k-1}}
~-~
(1+\frac{1}{p})\alpha_p.
\end{eqnarray*}
Since for any prime $p \ge 11$, one has $\alpha_p < \frac{6}{p}$ and hence
$$
\frac{\lambda_F(p^m)}{p^{m(k-1)}} 
~\geq~
1 - \frac{1}{p^{1/2}}(2 + \frac{5}{p^{1/2}} + \frac{2}{p} + \frac{6}{p^{3/2}}).
$$
Thus except for finitely many primes $p$, 
there exists an absolute constant $c_{10}>0$ such that 
for all $m\in \N$,
\begin{eqnarray}\label{eq12}
\frac{\lambda_F(p^m)}{p^{m(k-1)}} 
& \geq & 1-\frac{c_{10}}{p^{1/2}}
\phantom{m}\text{with   }\frac{c_{10}}{p^{1/2}}<1.
\end{eqnarray}
It is easy to see that one can choose $c_{10} = 3.5$ so that
the inequality \eqref{eq12} happens for any prime~$p \ge~17$.
Let 
$$
T :=\{p \in \mathcal{P} : \text{ the inequality }
\eqref{eq12}\text{ holds}\}
$$ 
and $n\in\N$ be any natural number whose prime
divisors are in $T$. As before, writing 
$$
n = \prod_{1 \le i \le t}p_i^{m_i}
$$ 
with $m_i >0$ and $p_1 < \cdots < p_t$, 
we have 
\begin{eqnarray*}
\frac{\lambda_F(n)}{n^{k-1}} 
~~\geq~
\prod_{1\leq i\leq t}(1 - \frac{c_{10}}{p_i^{1/2}})
&=& 
\text{exp}\left(\sum_{1\leq i\leq t}
\log(1 - \frac{c_{10}}{p_i^{1/2}})\right) \\
& \geq & 
\text{exp}\left(-c_{11}\sum_{1\leq i\leq t}
\frac{1}{p_i^{1/2}}\right)\\
& \geq & 
\text{exp}\left(-c_{11}\sum_{1\leq i\leq t}
\frac{1}{i^{1/2}}\right)\geq\text{exp}\left(-c_{12}~ t^{1/2}
\right),
\end{eqnarray*}
where $c_{11}, ~c_{12} >0$ are absolute constants.
Again since $t = \nu(n) \ll \frac{\log n}{\log\log n}$
for $n \gg 1$, we have for such $n\in \N$ with $n\ge 3$,
\begin{equation}\label{eq13}
\frac{\lambda_F(n)}{n^{k-1}}\geq 
\text{exp}\left(-c_3\sqrt{\frac{\log n}{\log\log n}}\right),
\end{equation}
where $c_3>0$ is an absolute constant.
Note that \eqref{eq13} holds if all the prime divisors 
of $n$ are in the set $T$. Now if $n \in \N$ is such that
$p|n\Rightarrow p\notin T$, then we use Hecke relation
$$
a(p^{n+1}) = a(p)a(p^n) - p^{2k-3}a(p^{n-1})
$$
for $n\in \N$ and explicit calculations using Mathematica.
In particular, we show that
\begin{equation}\label{eq14}
\frac{\lambda_F(n)}{n^{k-1}}
~\geq~ c_2,
\end{equation}
where $c_2 > 0$ is an explicit constant.
Combining \eqref{eq13} and \eqref{eq14}, we now get
$$
\lambda_F(n)
\geq
 c_2 n^{k-1}\text{exp}
\left(-c_3\sqrt{\frac{\log n}{\log\log n}}\right)
$$
for any natural number $n \in \N$ with $n\ge 3$.

\subsection{Proof of Corollary \ref{co}}
Since $\lambda_F$ is a non-zero multiplicative
function (see \cite{AA} and \cite{TO}), we have 
$\lambda_F(1) =1 > 0$. Also we know
that
$$
\lambda_F(2) \ge \frac{3}{2} - \sqrt{2} > 0.
$$
Now by applying Theorem \ref{th3},
we have our corollary.

\section{Proof of theorem \ref{th4}}\label{sec6}

Notations are as in the previous sections. For any
prime $p$, we get 
$$
\lambda_F(p)=p^{k-1}(1+\frac{1}{p}+\frac{a(p)}{p^{k-1}}).
$$
Write $b(p) = a(p) / p^{k-3/2}$. For any absolute constant
$0 < \tilde{\beta} < 2$,  consider the sums
\begin{eqnarray*}
&& 
S(x) := \sum_{p \le x} (b(p) + \tilde{\beta}) ~(b(p) -2)
\phantom{m}\text{and}\phantom{m}
S^+(x) := \sum_{p \le x, \atop b(p) < -\tilde{\beta}}
(b(p) + \tilde{\beta}) ~(b(p) -2). 
\end{eqnarray*}
Note that 
$$
S(x) \le S^+(x) \le 16 ~\#\{ p \in \mathcal{P} : b(p) < - \tilde{\beta} \} ~\log x.
$$
Then using the estimates (see pages 43 and 135 of \cite{IK} and 
Theorem 2 of  \cite{RR})  
\begin{eqnarray*}
&&
\sum_{p \le x} b(p) \log p  \ll  x~{\rm exp}(-\kappa \sqrt{\log x}), 
\phantom{m}
\sum_{p \le x} b^2(p) \log p  \sim x 
\phantom{m}
\text{and}
\phantom{m}
\sum_{ p\le x} \log p \ll x,
\end{eqnarray*}
where $\kappa >0$ is an absolute constant
and proceeding along the lines of the proof of lemma \ref{lem6}, 
one can show there exists an absolute 
constant  $0< \beta_1< 2$ such that 
$$
B:=\{p: a(p)< - \beta_1\cdot p^{k-3/2}\}
$$
has positive lower density. Let us take 
$$
n_x =\prod_{x < p \leq 2x \atop p\in B}p,
$$ 
where $x$ is sufficiently large so that 
$\frac{2}{{\sqrt{x}}} < \beta_1$. Then we have 
\begin{eqnarray*}
\frac{\lambda_F(n_x)}{n_x^{k-1}}
=  
\prod_{x < p \leq 2x,\atop p\in B}
\left(1+\frac{1}{p} +\frac{a(p)}{p^{k-1}}\right)
&\leq& 
\prod_{x < p\leq 2x,\atop p\in B}
\left(1+\frac{1}{p}+ \frac{-\beta_1}{p^{1/2}}\right) \\
& \leq & 
\text{exp}\left[\sum_{x < p\leq 2x,\atop p\in B}
\log(1-\frac{\beta_1}{2p^{1/2}})\right]\\
& \leq & 
\text{exp}\left(-c_{13}\sum_{x < p \leq 2x,\atop p\in B}
\frac{1}{p^{1/2}}\right),
\end{eqnarray*}
where $c_{13} >0$ is an absolute constant.
Since the set $B$ has positive lower density, as in 
section~\ref{sec3}, we get 
\begin{eqnarray*}
\frac{\lambda_F(n_x)}{n_x^{k-1}}
\leq
\text{exp}\left(-c_{15}\frac{\sqrt{x}}{\log x}\right)
\leq
\text{exp}\left(-c_4\frac{\sqrt{\log n_x}}{\log\log n_x}\right),
\end{eqnarray*}
where $c_{15} >0$ is an absolute constant.
Thus for given any natural number $M$, there exists 
a natural number $n$ with $n>M$ such that 
$$
\frac{\lambda_F(n)}{n^{k-1}}\leq \text{exp}
\left(-c_4\frac{\sqrt{\log n}}{\log\log n}\right).
$$
Hence we have the result.

\section{Proof of theorem \ref{th5}}

Recall that 
for any $m\in\N$ and any prime $p$, one has (see equation
\eqref{imp})
\begin{equation*}
\frac{\lambda_F(p^m)}{p^{m(k-1)}}
=
1+\frac{1}{p}
~+~
(1+\frac{1}{p})
\sum_{\ell=1}^{m-1}\frac{a(p^\ell)}{p^{\ell(k-1)}}
+
\frac{a(p^m)}{p^{m(k-1)}}
\end{equation*}
with the convention that an empty sum is
zero. Note that the series 
$$
\sum_{\ell=1}^\infty \frac{a(p^\ell)}{p^{\ell(k-1)}}
$$
is absolutely convergent (see section \ref{sec4} for details). 
This implies 
that the sequence 
$$
\left\{\frac{\lambda_F(p^m)}{p^{m(k-1)}}\right\}_{m\in\N}
$$ 
is convergent. Further there exist absolute constants $e_1, e_2>0$ 
such that 
\begin{eqnarray}\label{eq15}
1+\frac{e_1}{p^{1/2}}
\leq 
\frac{\lambda_F(p^m)}{p^{m(k-1)}} 
\leq  
1+\frac{e_2}{p^{1/2}}
\end{eqnarray}
holds for all but finitely many primes $p\in A$. 
Indeed, the upper bound is a consequence of section \ref{sec4}
whereas the lower bound follows from the fact that
$\alpha_p \le 6/p$ for $p \ge 11$ and 
primes $p \in A$ has the property
that $a(p) > \beta. ~p^{k-3/2}$ with absolute
constant $\beta$ (see section \ref{sec3}, section \ref{sec4} and 
section~\ref{sec5}).

Let us choose a prime $p_1\in A$ such that \eqref{eq15} holds. 
Since \eqref{eq15} is true for all but finitely many $p\in A$, 
we can choose $p_2\in A$ such that $p_2>p_1$ and
\begin{eqnarray*}
1+\frac{e_1}{p_2^{1/2}}
~\leq~~ 
\frac{\lambda_F(p_2^m)}{p_2^{m(k-1)}} 
~\leq~~  
1+\frac{e_2}{p_2^{1/2}}
<  
1+ \frac{e_1}{2p_1^{1/2}}.
\end{eqnarray*}
Proceeding in this way, we get a sequence 
$\{p_n\}_{n\in\N}$ such that each
$$
\lim_{m\rightarrow\infty}
\frac{\lambda_F(p_n^m)}{p_n^{m(k-1)}}
>1
\phantom{m}
\text{and} 
\phantom{m}
\lim_{m\rightarrow\infty}
\frac{\lambda_F(p_i^m)}{p_i^{m(k-1)}}
\ne
\lim_{m\rightarrow\infty}
\frac{\lambda_F(p_j^m)}{p_j^{m(k-1)}}
$$
for any $i\ne j$. Thus
there are infinitely many limit 
points of the sequence 
$\{\frac{\lambda_F(n)}{n^{k-1}}\}_{n\in\N}$ 
which are~$>1$.

Considering the set $B$ (see section \ref{sec6}) and
arguing as above, we can show there is
a sequence 
$\{p_n\}_{n\in\N} \subset B$ for which
$$
\lim_{m\rightarrow\infty}
\frac{\lambda_F(p_n^m)}{p_n^{m(k-1)}}
<1
\phantom{m}
\text{and} 
\phantom{m}
\lim_{m\rightarrow\infty}
\frac{\lambda_F(p_i^m)}{p_i^{m(k-1)}}
\ne
\lim_{m\rightarrow\infty}
\frac{\lambda_F(p_j^m)}{p_j^{m(k-1)}}
$$
for any $i\ne j$. This completes the proof.

\bigskip
\noindent
{\bf Acknowledgments:} The third author would like to thank 
the Institute of Mathematical Sciences for providing excellent 
working atmosphere where the work was done.
The authors would like to thank Purusottam Rath for asking
the possibility of infinitude of limit points of the
sequence considered in Theorem \ref{th5} and going
through an earlier version of the paper.

\end{document}